\documentclass{amsart}
\usepackage{amssymb,latexsym,enumerate,amscd}

\swapnumbers

\makeatletter
\newcommand{\shortto}{\mathchoice{\scriptstyle\to}{\scriptstyle\to}
{\scriptscriptstyle\to}{\scriptscriptstyle\to}}
\newcommand{\overto}[2]{\lower.7pt\vbox{\baselineskip=0pt\lineskip-.1pt
\ialign{$\m@th#1\hfil##\hfil$\crcr#2\crcr\shortto\crcr}}}
\newcommand{\largesim}{\mathchoice{\sim}{\sim}{\sim}{\scriptstyle\sim}}
\newcommand{\isoto}{\mathrel{\mathpalette\overto\largesim}}
\makeatother

\newtheorem{Thm}[equation]{Theorem}
\newtheorem{Prop}[equation]{Proposition}
\newtheorem{Lem}[equation]{Lemma}
\newtheorem{Cor}[equation]{Corollary}

\theoremstyle{remark}

\newtheorem*{Rem*}{Remark}
\theoremstyle{definition}

\newtheorem*{Not*}{Notation}
\newtheorem{Def}[equation]{Definition}
\numberwithin{equation}{section}
\DeclareMathOperator{\pr}{pr}
\DeclareMathOperator{\id}{id}
\DeclareMathOperator{\gap}{gap}
\DeclareMathOperator{\adef}{adef}
\DeclareMathOperator{\End}{End}

\begin{document}

\date{%
Mon Dec 24 16:37:04 EST 2007}

\title[Generation Gaps and Free Products]
{Generation Gaps and Abelianised Defects of Free Products}

\author[K. W. Gruenberg]{Karl W. Gruenberg}
\address{School of Mathematical Sciences \\
Queen Mary, University of London \\
Mile End Road \\
London E1 4NS \\
England}
\author[P. A. Linnell]{Peter A. Linnell}
\address{Department of Mathematics \\
Virginia Tech \\
Blacksburg \\
VA 24061-0123\\
USA}
\email{linnell@math.vt.edu}
\urladdr{http://www.math.vt.edu/people/linnell/}
\thanks{The first author died on 10 October 2007.
The second author gratefully acknowledges financial support
from the EPSRC}

\begin{abstract}
Let $G$ be a group of the form $G_1 * \dots *G_n$, the free product
of $n$ subgroups, and let $M$ be a $\mathbb{Z}G$-module of the form
$\bigoplus_{i=1}^n M_i \otimes_{\mathbb{Z}G_i} \mathbb{Z}G$.  We
shall give formulae in various situations
for $d_{\mathbb{Z}G}(M)$, the minimum number of
elements required to generate $M$.  In particular if $C_1,C_2$ are
nontrivial finite cyclic groups of coprime orders,
$G = (C_1 \times \mathbb{Z}) * (C_2 \times \mathbb{Z})$
and $F/R \cong G$ is the free presentation obtained from
the natural free presentations of the two factors, then the number
of generators of the relation module, $d_{\mathbb{Z}G}(R/R')$, is
three.  It seems plausible that
the minimum number of relators of $G$ should be 4, and this
would give a finitely presented group with positive relation gap.
However we cannot prove this last statement.
\end{abstract}



\maketitle

\section{Introduction} \label{Sintroduction}

Given a group $G$ and a $\mathbb{Z}G$-module $M$,
let $d(G)$ denote the minimal number of generators
of $G$, let $d_G(M)$ denote the minimal number of
$\mathbb{Z}G$-module generators of $M$,
and let $\Delta G$ denote the augmentation ideal of $G$, that is the
ideal of $\mathbb{Z}G$ generated by $\{g-1 \mid g \in G\}$.
Then the \emph{generation gap} of the finitely generated group $G$ is
the difference
\[
\gap G := d(G) - d_G(\Delta G).
\]
If $G$ is finitely presentable and $R\hookrightarrow F
\twoheadrightarrow G$ is a finite free presentation (meaning
$F$ is a finitely generated free group and $R$ has a finite
number of generators as a normal subgroup of $F$), then the
\emph{abelianised defect} $\adef(F/R)$ of the free presentation
is the difference $d_G(R/R') - d(F)$, where
$R'$ is the commutator group of $R$ and the free abelian group
$R/R'$ is viewed as a $\mathbb{Z}G$-module (the relation module of
the presentation).  A fair amount is known about the numerical
functions gap and adef for finite groups.  Our purpose here is to
make a start with their study for infinite groups.  In particular, to
indicate how results on finite groups can be used to discuss free
products.  In the next section, after reviewing some relevant facts of
the finite theory, we state and discuss our main results.  We mention
just two striking consequences:
\begin{Prop} \label{Pfinitegroups}
If $G$ is the free product $G_1 * \dots * G_n$ of the finite
groups $G_1, \dots, G_n$ of mutually coprime orders, then $\gap(G) =
0$ if and only if there exists a positive integer $k$ such that
$\gap (G_k) = 0$ and $G_i$ is cyclic for all $i \ne k$.
\end{Prop}

\begin{Prop} \label{Pinfiniterel}
Let $G_1,G_2$ be nontrivial finite
cyclic groups of coprime orders, let $C_1,C_2$ be infinite cyclic
groups, and let $G = (G_1 \times C_1) * (G_2 \times C_2)$.  If $F/R
\isoto G$ is the free presentation obtained from the
natural $2$-generator free presentations
of the two factors, then $d_G(R/R')=3$.
\end{Prop}
The second of these results could have an important consequence for
the relation gap problem: given a finite free presentation
$R\hookrightarrow F \twoheadrightarrow G$
of a group $G$, then the \emph{relation gap} of the presentation is
defined to be $d_F(R) - d_G(R/R')$, where $d_F(R)$ denotes the
minimum number of elements needed
to generate $R$ as a normal subgroup of $F$.
No examples are known of positive relation gap.  The group of
Proposition \ref{Pinfiniterel} above needs at most
4 (obvious) relations.  If 4
really equals $d_F(R)$, then $F/R$ would have positive relation gap.
Proposition \ref{Pinfiniterel} was mentioned in the interesting
survey article \cite{Harlander00} by Jens Harlander.

A similar result was proved in \cite[Proposition
1]{BridsonTweedale07}, where consequences to the ``$D(2)$
conjecture" \cite{Johnson03} were discussed.  We shall show how to
obtain \cite[Propositions 1,2]{BridsonTweedale07} from Theorem
\ref{Tmain2} at the end of Section \ref{SIFF}.

The first author Karl Gruenberg died on 10 October 2007.  At that
time this paper was more or less complete; mistakes that remain are
entirely the fault of the second author, Peter Linnell.  Karl was a
close and generous friend to not only the second author, but many
other mathematicians.  His influence through his research and
teaching will be greatly missed.

\section{Review and Overview} \label{Sreview}

Among finite groups, generation gap zero is the rule rather than the
exception: all simple groups (indeed all 2-generator groups) and all
soluble groups have generation gap zero.  However there do exist
groups of arbitrarily large generation gap.  The smallest group with
nonzero generation gap is known to be the direct product of 20 copies
of the alternating group of degree 5 \cite{VoltaLucchini99}.
For more on these matters, see \cite[\S 13]{Gruenberg79}.

Finite groups also behave well with respect to abelianised defect: on
a given finite group $G$, the function $\adef(F/R)$ as $F/R$ varies
over all finite free presentations of $G$ is constant
(cf.~\cite[\S 14]{Gruenberg79}).  Hence $\adef(G)$ is an unambiguous
notation for this integer.  For finitely presentable groups $G$, the
set of $\adef(F/R)$ for varying finite free presentations $F/R$ is
bounded below; so the lower bound may be written $\adef(G)$.
A free presentation $F/R \isoto G$ is called
\emph{minimal} if $d(F) = d(G)$.

If $G = G_1* \dots *G_n$ and we are given free presentations
$F_i/R_i \isoto G_i$ for $i= 1,\dots,n$, then $F = *F_i$ projects
onto $G$ with kernel the normal closure of $R_1,\dots,R_n$.

If $R$ is a relation group, we shall always
abbreviate the resulting relation module $R/R'$ as
$\overline{R}$.  When the relation group is $R_i$, we shall use the
streamlined notation $\overline{R}_i$ for $\overline{R_i}$.

\begin{Thm} \label{Tfinitegroups}
Let $G= G_1* \dots *G_n$ with all the $G_i$ finite groups
of mutually coprime orders.  Then
\begin{enumerate} [\normalfont(a)]
\item
\[
d_G(\Delta G) = \max_k\{d_{G_k}(\Delta G_k)\} + n - 1
\]
whence
\[
\gap(G) = \min_k\{\gap (G_k) + \sum_{i\ne k} (d(G_i)-1)\}.
\]

\item
Let $F/R \isoto G$ be the free presentation constructed
from free presentations $F_i/R_i \isoto G_i$ of
the free factors.  Then
\[
d_G(\overline{R}) = \max_k\{d_{G_k}(\overline{R}_k) +
\sum_{i\ne k} d(F_i)\},
\]
whence $\adef(F/R) = \max_k\{\adef(G_k)\}$.
\end{enumerate}
\end{Thm}

An immediate consequence of Theorem \ref{Tfinitegroups}(a)
is Proposition \ref{Pfinitegroups}.  Theorem
\ref{Tfinitegroups}(b) tells us, for example,
that if each $G_i$ is an elementary abelian $p_i$-group of
rank 2 and $d(F_i) =2$, then $d_G(\overline{R}) =3 + 2(n-1)
= 2n+1$ because here
(clearly) $d_{G_i}(\overline{R}_i) = 3$.  One might expect $d_F(R)$
to be $3n$ since $R$ is the normal closure of $R_1,\dots,R_n$ and
$d_{F_i}(R_i) = 3$ for all $i$.  However in the case $n=2$,
Hog-Angeloni, Lustig and Metzler
showed that $d_F(R) = 5$ \cite[Addendum to Theorem 3 on
p.~163]{HLM85} and hence the relation gap of $F/R$ is zero.
Their results yield for general $n$ that $d_F(R) \le 3n-1$.  If
there exists $n>2$ such that $d_F(R) = 3n-1$, then the group involved
will have positive relation gap, namely $3n-1 - (2n+1) =
n-2 > 0$.

Theorem \ref{Tfinitegroups} follows from the general
Theorem \ref{Tmain2} stated and proved in \S\ref{Smain}.
All our results on gap and adef are consequences of this.
Theorem \ref{Tmain2} is a join of two results, the first of
which, Theorem \ref{Tmain1}, is our major contribution:
it gives an upper bound for the minimal number of generators
of a $\mathbb Z G$-module $M$, where $G$ has subgroups $G_1,
\dots, G_n$ and $M$ contains ``\emph{good}" $\mathbb{Z}
G_i$-submodules $M_i$ that together $G$-generate $M$.
When $G= G_1* \dots *G_n$ and $M = \bigoplus_{i=1}^n
(M_i \otimes_{\mathbb Z G_i} \mathbb Z G)$, our upper bound
for $d_G(M)$ becomes a lower bound by a remarkable theorem
of George Bergman \cite{Bergman74}.  The restriction to what we
call ``good" modules (cf.\ Definition \ref{Dgood}) is probably
unavoidable.  Fortunately it is no restriction at all for finite
groups.  The proof of this has no bearing on any of our other
arguments and we therefore postpone it to an appendix.

The following two theorems look somewhat like Theorem
\ref{Tfinitegroups}, but here the factors are infinite groups and
this leads to difficulties with the ``good" property of modules.
\begin{Thm} \label{Tinfinitegap}
Let $H = H_1 * \dots * H_n$, where each $H_i = G_i \times A_i$,
$G_i$ is a finite nilpotent group, $A_i$ is a finitely generated
nilpotent group with $A/A'$ torsion free and the
orders of $G_1, \dots, G_n$ are mutually coprime.  Then
\begin{enumerate}[\normalfont(a)]
\item
\[
d_H(\Delta H) = \max_k\{d(H_k) + \sum_{i\ne k}(d(A_i) +
\delta_{A_i,1})\},
\]
where $\delta_{A_i,1} = 1$ when $A_i = 1$ (equivalently $d(A_i) = 0$)
and $\delta_{A_i,1} = 0$, when $A_i \ne 1$.

\item
$\gap(H) = \min_k\{\sum_{i\ne k} (d(G_i) - \delta_{A_i,1})\}$.
\end{enumerate}
\end{Thm}
We conclude that a
\emph{free product of finitely generated abelian groups, with
torsion groups of mutually coprime orders,
has generation gap zero if, and only if, every factor but one is
torsion free or finite cyclic}.  Theorem
\ref{Tinfinitegap} is an immediate corollary of the more general
Theorem \ref{Taugmentation} in \S\ref{SIFF}.

Finally, we state the theorem of which Proposition
\ref{Pinfiniterel} is a special case.
\begin{Thm} \label{Tgeneralrelation}
Let $H = H_1 * \dots * H_n$, where each $H_i = G_i \times C_i$ with
$G_i$ finite cyclic and $C_i$ infinite cyclic.  Assume that the
orders of $G_1,\dots,G_n$ are mutually coprime.  Take natural
minimal free presentations $F_i/R_i \isoto H_i$ and let $F/R
\isoto H$ be the resulting free presentation of $H$.
Then $d_H(\overline{R}) = n+1$.
\end{Thm}
Theorem \ref{Tgeneralrelation} is an immediate corollary of
Theorem \ref{Tfinal}.  It seems
probable to us that Theorem \ref{Tgeneralrelation} remains true (in
an obviously modified enunciation) when each $C_i$ is allowed to be
free abelian of finite rank and $G_i$ is a suitably restricted finite
group.  Just what may be needed to establish such a generalisation
will become clear in \S\ref{SIFF}.

\section{Preliminaries} \label{SPreliminaries}

We explain some notation and collect a number of results, mostly
well-known, that are
needed later.  Modules will be right modules and mappings will be
written on the left.
If $V$ is a finitely generated module over a ring $S$,
we write $d_S(V)$ for the minimum number of elements needed to
generate $V$ as an $S$-module.  Of course if $X \subseteq V$, then
$XS$ indicates the $S$-submodule of $V$ generated by $X$.  In the
case $S=kG$, the group ring of the
group $G$ over the commutative ring $k$, we usually write $d_G(V)$
instead of $d_{kG}(V)$.  If $G$ is finite, we let $\hat{G} =
\sum_{g \in G} g \in kG$.
If $K$ is a commutative ring containing $k$,
then $KV$ means $K \otimes_k V$, so if $V$ is a $\mathbb{Z}G$-module,
$d_G(\mathbb{Q}V)$ means $d_{\mathbb{Q}G}(\mathbb{Q}
\otimes_{\mathbb{Z}} V)$.  For a prime $p$, we shall let
$\mathbb{F}_p$ indicate the field with $p$ elements.
The commutator subgroup of $G$
will be denoted by $G'$ and we let $G^n = \langle g^n \mid g \in
G\rangle$ for any positive integer $n$.
If $H$ is a group, then $\delta_{G,H}
\in \mathbb{Z}$ is the Kronecker $\delta$:
so $\delta_{G,H} = 1$ if
$G \cong H$ and $\delta_{G,H} = 0$ if $G \ncong H$.
Thus the $\delta$ in Theorem 2.2 is the Kronecker delta.

If $\pi$ is a set of primes, then $\pi'$ denotes the complementary
set.  Let $\pi(G)$ denote the set of primes dividing the orders of
the finite subgroups of $G$.  If $G$ is torsion free, set
$\pi(G) = \emptyset$.  Let $M$ be a finitely generated
$\mathbb{Z}G$-module.  We call $M$ a \emph{Swan module} if there
exists a prime $p$ such that $d_G(M) = d_G(M/pM)$; in that case
we call $p$ an \emph{$M$-prime}.  If $G$ is finite and $M$ is a
$\mathbb{Z}G$-lattice, then $M$-primes, if they exist, are divisors
of $|G|$ (cf.~\ref{Lmodulegenerators}).  There is a well developed
theory of Swan $\mathbb{Z}G$-lattices when $G$ is finite:
cf.~\cite[\S 7.1]{Gruenberg76} for the basics.  Augmentation ideals
and relation modules of finite groups are Swan lattices.  A
\emph{relation prime} of $G$ is an $\overline{R}$-prime, where
$\overline{R}$ is the relation module coming from a
finite free presentation of $G$.  When $G$ is finite, this
notion is independent of the chosen presentation: if
$\overline{R}_1,\overline{R}_2$ are relation modules,
then the $\overline{R}_1$-primes coincide with the
$\overline{R}_2$-primes: cf.~\cite[Lemma 7.13]{Gruenberg76}.

\begin{Lem} \label{Lgap0}
Let $H$ be a finitely generated group and let $G$ be a homomorphic
image of $H$ such that $d(G) = d(H)$.  If $\gap(G) = 0$ and
$\Delta G$ is a Swan module for $G$, then $\gap(H) = 0$ and
$\Delta H$ is a Swan module for $H$.  Moreover every $\Delta G$-prime
is a $\Delta H$-prime.
\end{Lem}
\begin{proof}
If $p$ is a $\Delta G$-prime, then $d(H) = d(G) =
d_G(\Delta G) = d_G(\Delta G/p\Delta G) \le d_H(\Delta H/p\Delta H)
\le d_H(\Delta H) \le d(H)$.
\end{proof}

\begin{Cor} \label{Cgap0}
If $H$ is a finitely generated nilpotent group, then $\gap(H) = 0$
and $\Delta H$ is a Swan module.  When $H/H'$ is torsion free,
every prime is a $\Delta H$-prime; otherwise every $\Delta H$-prime
belongs to $\pi(H/H')$.
\end{Cor}
\begin{proof}
There exists a prime $p$ such that $G := H/H'H^p$ satisfies $d(G)
= d(H)$.  Since $G$ is an elementary abelian $p$-group, it has
generation gap zero.  When $H/H'$ is torsion free, any prime will
do for $p$; but if $H/H'$ has torsion subgroup $T \ne 1$, then $p$
must divide $|T|$.
\end{proof}

\begin{Lem} \label{Lmodulegenerators}
Let $G$ be a finite group and let $M$ be a finitely generated
$\mathbb{Z}G$-module.  Then $d_G(\mathbb{Q}M) \leq d_G(M/pM)$ for
every prime $p$, while $d_G(\mathbb{Q}M) = d_G(M/pM)$ for every
$p \notin \pi(G) \cup \pi(M)$.
\end{Lem}
\begin{proof}
First we show that $d_G(\mathbb{Q}M) \le d_G(M/Mp)$.  Choose $x_1,
\dots, x_n \in M$ to generate $M$ as a $\mathbb{Z}G$-module modulo
$p$.  Then the elements $x_ig$ for $1 \le i \le n$ and $g \in G$
generate $M$ as a $\mathbb{Z}$-module modulo $p$.  Since $M$ is
finitely generated as a $\mathbb{Z}$-module, we see that the
elements $x_ig$ generate $\mathbb{Q}M$ as a $\mathbb{Q}$-module.
Hence the elements $x_i$ generate $\mathbb{Q}M$ as a
$\mathbb{Q}G$-module.

Next we show that $d_G(M/Mp) \le d_G(\mathbb{Q}M)$ when
$p \notin \pi(G) \cup \pi(M)$; this will complete
the proof.  Choose $x_1,\dots,x_n \in M$ to generate $\mathbb{Q}M$ as
a $\mathbb{Q}G$-module.  Write $N = x_1\mathbb{Z}G + \dots +
x_n\mathbb{Z}G$.  Then $\mathbb{Q}N = \mathbb{Q}M$, consequently
$\mathbb{C} \otimes_{\mathbb{Z}} N \cong \mathbb{C}
\otimes_{\mathbb{Z}} M$.  If $M'$ is the torsion group of $M$ and
$\overline{M} = M/M'$, then $\overline{M}$ is a
$\mathbb{Z}G$-lattice and $\overline{M}/p\overline{M} \cong
M/pM$ since $p$ is prime to $|M'|$.  Let $k$ denote the algebraic
closure of the field $\mathbb{F}_p$.  By Brauer character
theory, we see by looking at the decomposition numbers that
$k \otimes_{\mathbb{Z}} M$ and $k\otimes_{\mathbb{Z}} N$ have the same
composition factors as
$kG$-modules.  Since $p$ is prime to $|G|$, all $kG$-modules are
completely reducible and we deduce that $k\otimes_{\mathbb{Z}} M
\cong k \otimes_{\mathbb{Z}} N$.  Therefore $M/Mp \cong N/Np$ as
$\mathbb{F}_pG$-modules by the Noether-Deuring theorem.
Thus $d_G(M/pM) \le n$ as required.
\end{proof}

\section{The Main Technical Result} \label{Smain}

Given a finitely generated $\mathbb{Z}G$-module $M$ and a finite set
of primes $\pi$, let $(X_p)_{p\in \pi}$ be a family of subsets of $M$
such that $X_p$ is a minimal $\mathbb{Z}G$-generating set of $M$
modulo $pM$ and $X_p \supseteq X_q$ if, and only if, $d_G(M/pM)
\geq d_G(M/qM)$.  We call $(X_p)_{p\in \pi}$ a
\emph{nested $\pi$-family in} $M$.  Note that the nested structure
ensures that if $X_q$ is the minimal element of $(X_p)_{p\in \pi}$,
then for every $p \in \pi$,
\begin{equation} \label{Enest}
d_G\bigl(M/(pM + X_q\mathbb Z G)\bigr) = d_G(M/pM) - d_G(M/qM).
\end{equation}
If $\pi = \{p_1,\dots, p_t\}$, we will write $X_i = X_{p_i}$ for $1
\le i \le t$ and assume that $X_1,\subseteq \dots,\subseteq X_t$.

\begin{Lem} \label{Lnest}
For every $M$ and $\pi$ there exists a nested $\pi$-family in $M$.
\end{Lem}
\begin{proof}
Let $\pi = \{ p_1, \dots, p_t\}$ be ordered so that $d_G(M/p_iM)
\geq d_G(M/p_jM)$ if, and only if, $i \leq j$.  Choose a subset
$Y_i$ of $M$ so that its image in $M/p_iM$ is a minimal
$\mathbb{Z}G$-generating set of $M/p_iM$.  Then $|Y_1| \geq |Y_2|
\geq \dots \geq |Y_t|$.

Let $d = |Y_1|$ and enumerate the elements in each $Y_i$ using
zeros if necessary at the end: $Y_i = \{y_{i,1}, y_{i,2},
\dots, y_{i,d}\}$.  Let $p_i' = \prod_{j\ne i}
p_j$.  Then $(p_i, p_i') = 1$ and $p_im_i + p_i'm_i' = 1$ for
suitable integers $m_i, m_i'$.  We define
\[
x_k = \sum_{i=1}^t p_i'm_i' y_{i,k} \qquad (1\leq k \leq d).
\]
Thus $x_k \equiv p_j'm_j'y_{j,k} \equiv y_{j,k}$
modulo $p_jM$, whence $X_i = \{ x_1, \dots, x_{d_i}\}$
with $d_i = |Y_i|$ gives the required nested $\pi$-family
$X_1 \supseteq X_2 \supseteq \dots \supseteq X_t$.
\end{proof}

\begin{Def} \label{Dgood}
Let $G$ be a group and let $M$ be a finitely generated
$\mathbb{Z}G$-module.  We say $M$ is \emph{good} for $G$ if there
exists a finite set of primes $\pi(G,M)$ such that
\begin{enumerate}[\normalfont(g1)]
\item
$d_G(M/pM)$ has constant value $\delta_G(M)$ for all $p
\notin \pi(G,M)$, while $\delta_G(M) \leq d_G(M/pM)$
for all $p \in \pi(G,M)$;
\item
for every finite set of primes $\pi$ properly containing
$\pi(G,M)$ there exists a nested $\pi$-family whose minimal
element $X$ satisfies
\[
d_G\bigl(M/(pM + X\mathbb Z G)\bigr)\leq 1\quad \text{ for all }
p \notin \pi,
\]
and $M/X\mathbb{Z}G$ is a torsion group of finite exponent.
\end{enumerate}
\end{Def}

It is then automatic that $|X| = \delta_G(M)$.
Also if $p\in \pi'$, then $d_G\bigl(M/(pM + X \mathbb{Z}G) \bigr)
= 1$ if, and only if, $ p \in \pi(M/X \mathbb Z G)$ (because $M\ne
pM + X \mathbb{Z}G$ if, and only if, $p \in \pi(M/X\mathbb{Z}G)$).

\begin{Thm} \label{Tmain1}
Let $G$ be a group containing subgroups
$G_1,\dots,G_n$, let $M_i$ be a good $\mathbb{Z}G_i$-module and let
$M$ be a $\mathbb{Z}G$-module containing $M_1,\dots, M_n$ and
generated by them as a $\mathbb{Z}G$-module.  Assume that either all
the $M_i$ are Swan modules or that $d_{G_i} (M_i/pM_i) > \delta _{G_i}
(M_i)$ for some integer $i$ and some prime $p$.  Then
\[
d_G(M) \leq \max_p\big\{ \sum_{i=1}^n d_{G_i}(M_i/pM_i)\big \}.
\]

\end{Thm}
\begin{proof}
Set $\delta = \max_p \{ \sum_{i=1}^n d_{G_i}(M_i/pM_i) \}$,
$e_i = \delta_{G_i}(M_i)$ and let $\pi$ be a finite set of primes
strictly containing $\pi(G_1, M_1) \cup \dots \cup \pi(G_n, M_n)$.
Choose a nested $\pi$-family in $M_1$ as required for (g2) and
let $X_1$ be the minimal element.  Set $N_1 = M_1/X_1\mathbb Z G_1$.
Now $\pi(N_1)$ is finite and $d_{G_1} (N_1/pN_1) \ne 0$ if, and only
if, $p \in \pi(N_1)$.  Define $\pi_2 = \pi \cup \pi(N_1)$ and choose
a nested $\pi_2$-family in $M_2$ as required for (g2).
Let its minimal element be $X_2$ and set $N_2 =
M_2/X_2\mathbb{Z}G_2$.  Continue in this manner, producing
$X_3, \dots, X_n$, with $N_i = M_i/X_i\mathbb Z G_i$ and $\pi_i
= \pi \cup \pi(N_1)\cup \dots \cup \pi(N_{i-1})$ for $i=2, \dots, n$.
This has ensured that if $i < j$ and $p \in \pi'$, then
$d_{G_i}(N_i/pN_i)$ and $d_{G_j}(N_j/pN_j)$ cannot both be
non-zero.  For suppose $d_{G_i}(N_i/pN_i) \ne 0$.  Then $p \in
\pi(N_i)$, whence $p\in \pi_j $ and then, since $X_j$ is the
minimal element in a nested $\pi_j$-family, $d_{G_j}(N_j/pN_j)
= d_{G_j}(M_j/pM_j) - \delta_{G_j}(M_j)$ by \eqref{Enest},
where the right hand side is zero by (g1), because $p\notin
\pi(G_j, M_j)$.  If $d_{G_j}(N_j/pN_j) \ne 0$, then $p\notin
\pi_j$ (otherwise $d_{G_j}(N_j/pN_j) = 0$ as we have just seen)
and so $p \in \pi(N_i)'$, which ensures $d_{G_i}(N_i/pN_i) = 0$.
Thus we have
\begin{equation} \label{Ecyclicmodp}
\sum_{i=1}^n d_{G_i}(N_i/pN_i) \leq 1 \qquad \text{ for every }
p \in \pi'.
\end{equation}

Suppose $\delta = \sum_{i=1}^n e_i$.  Since $e_i \leq d_{G_i}
(M_i/pM_i)$ for every $p$ by (g1), we see that
$e_i = d_{G_i}(M_i/pM_i)$ for every $p$ and every $i$.
In this case our hypothesis makes every $M_i$ a Swan module, and so
$e_i = d_{G_i} (M_i)$ by using an $M_i$-prime $p$.  Thus $\delta =
\sum_{i=1}^n d_{G_i}(M_i)$, which gives the result.
Therefore we may assume that $\delta > \sum_{i=1}^n e_i$.

For each positive integer
$r \le \delta - \sum_{i=1}^n e_i$, let $\mathcal{P}_r$
denote the set of primes $p \in \pi$ for which $r \le
\sum_{i=1}^n d_{G_i} (N_i/pN_i)$.  Thus $\mathcal{P}_1 \ne
\emptyset$ (because $e_i < d_{G_i}(M_i/pM_i)$ for some $i$ and
some $p\in \pi$), $\mathcal P_1 \supseteq \mathcal P_2 \supseteq
\dots \supseteq \mathcal P_s$, where $s = \delta -
\sum_{i=1}^n e_i$ and $\mathcal P_{s+1} = \emptyset$.

Recall $N_i \ne p N_i$ for some $p\in \pi'$ if and only if
$\pi(N_i) \cap \pi' \ne \emptyset$.  Let $\pi_0(N_i) = \pi(N_i)
\cap \pi'$ and let $I_0$ be the set of all $i$ such that
$\pi_0(N_i) \ne \emptyset$.  By \eqref{Ecyclicmodp}, the sets
$\pi_0(N_i)$, for $i \in I_0$, are disjoint.  If $i \in I_0$,
we know $d_{G_i}(N_i/pN_i) = 1$ for every $p \in \pi_0(N_i)$
(property (g2)), whence we may choose a $p$-element $u_{i,0}(p)$
which $\mathbb{Z}G_i$-generates $N_i$ modulo $pN_i$.  If $i
\notin I_0$ set $u_{i,0}(p) = 0$ for all $p$.

We now define inductively $\mathbb{Z}G_i$-quotients $N_{i,k}$ of
$N_i$ for $1\leq k\leq s$, starting with $N_{i,1} = N_i$, for
$1\leq i \leq n$.

The set of primes $p\in \pi$ so that $N_i \ne p N_i$ for some $i$
is precisely $\mathcal{P}_1$.  Let $\pi_1(N_i) = \pi(N_i) \cap \pi$
and let $I_1$ be all $i$ so that $\pi_1(N_i) \not= \emptyset$.
With $i \in I_1$, choose a $p$-element $u_{i,1}(p) \in N_i$ so that
$d_{G_i} \bigl(N_i/(pN_i + u_{i,1}(p)\mathbb{Z} G_i)\bigr) < d_{G_i}
(N_i/pN_i)$.  If $\pi_1(N_i) =\emptyset$, put $u_{i,1}(p) =0$ for all
$p$.  Let $\pi(N_{i,1})$ be the union of the disjoint sets
$\pi_0(N_i)$, $\pi_1(N_i)$ and define
\[
K_{i,1} = \{u_{i,0}(p), \, u_{i,1}(p) \mid p
\in \pi(N_{i,1})\}\mathbb{Z}G_i.
\]
Set $N_{i,2} = N_i/K_{i,1}$.

If $N_{i,2} \ne pN_{i,2}$, then $p \in \mathcal P_2$ (because then
$N_i\ne pN_i$, $p \in \pi_1(N_i)$ and $d_{G_i}(N_i/pN_i) \ge 2$).
Thus $\pi(N_{i,2}) \subseteq \mathcal{P}_2$.  Let $I_2$ be all $i$
for which $\pi(N_{i,2}) \ne \emptyset$.  If $p \in \pi(N_{i,2})$,
choose a $p$-element $u_{i,2}(p)\in N_i$ so that $d_{G_i}
\bigl(N_{i,2}/(pN_{i,2} + \bar u_{i,2}(p)\mathbb Z G_i)\bigr)
< d_{G_i} (N_{i,2}/pN_{i,2})$, where $\bar u_{i,2}(p) =
u_{i,2}(p)+pN_{i,2}$.  If $\pi(N_{i,2}) = \emptyset$, set
$u_{i,2}(p) =0$ for all $p$.  Define
\[
K_{i,2} = K_{i,1}+ \{u_{i,2}(p) \mid p\in
\pi(N_{i,2})\}\mathbb{Z}G_i
\]
and $N_{i,3} = N_i/K_{i,2}$.  Continue until we reach $N_i = K_{i,r}$,
for some $r \leq s$ (because $\pi(N_{i,k}) \subseteq \mathcal P_k$
for all $k$).

By construction, for each $i$ and each $p$, the set of all non-zero
elements $u_{i,k}(p)$ for $k \geq 0$ maps to a minimal set of
$\mathbb{F}_pG_i$-generators of $N_i/pN_i$.  When $u_{i,0}(p) \ne 0$,
the set consists just of this one element because then $p \notin \pi$.
The number of non-zero $u_{i,k}(p)$ for all $k \geq 0$ and all $i$
is $\sum_{i=1}^n\,d_{G_i}(N_i/pN_i)$.  Thus
\[
\max_p|\{ u_{i,k}(p) \mid k\geq 0 \text{ and all } i\}| \leq s.
\]

Let $N = \bigoplus_{i=1}^n N_i\otimes_{\mathbb{Z}G_i} \mathbb{Z}G$.
Now $A = \langle u_{i,k}(p) \mid \text{all } i,k,p\rangle$ is a
\emph{finite} additive subgroup of $N$ and $N = A\mathbb{Z}G$
(because $N_i = K_{i,s}$ for all $i$).  The subset $\{ u_{i,k}(p)
\mid \text{all } i,k\}$ generates $A(p)$, the Sylow $p$-subgroup of
$A$ and $d(A) = \max_p d\bigl( A(p)\bigr)$.  Hence $d_G(N) \leq
\max_p d\bigl( A(p)\bigr) \leq s$.

Thus finally, $d_G(M) \leq d_G(N) + |X_1|+\cdots +|X_n| \leq s
+ \sum e_i = \delta$.
\end{proof}

All our results on free products are consequences of the following
theorem which is an immediate consequence of Theorem \ref{Tmain1}
and a theorem of Bergman \cite[Theorem 2.3]{Bergman74}:

\begin{Thm} \label{Tmain2}
Let $G = G_1* \dots * G_n$, let $M_i$ be a good
$\mathbb{Z}G_i$-module for $1 \le i\le n$, and let
$M = \bigoplus_{i=1}^n M_i\otimes_{\mathbb{Z}G_i} \mathbb{Z}G$.  If
either all the $M_i$ are Swan modules or $d_{G_i}(M_i/pM_i) >
\delta_{G_i}(M_i)$ for some $i$ and some prime $p$, then
\[
d_G(M) = \max_p \{ \sum_{i=1}^n d_{G_i}(M_i/pM_i) \}.
\]
\end{Thm}
\begin{proof}
Theorem \ref{Tmain1} gives the inequality $\leq$, while Bergman's
theorem tells us that, for each $p$, we have $\sum_{i=1}^n
d_{G_i}(M_i/pM_i) = d_G(M/pM)$.
\end{proof}
Note that, by Bergman's theorem, the conclusion of Theorem
\ref{Tmain2} is equivalent to the statement that
\emph{$M$ is a Swan $\mathbb{Z}G$-module}.

\section{Finite Free Factors} \label{SFFF}

\emph{If $G$ is a finite group, then every finitely generated
$\mathbb{Z}G$-module $M$ is good for $G$} if we take $\pi(G,M)
= \pi(G) \cup \pi(M)$.  We already know that (g1) holds: this is
Lemma \ref{Lmodulegenerators} with $\delta_G(M) = d_G(\mathbb Q M)$.
We postpone the proof of (g2) to the Appendix.

For the rest of this section we assume that
$G_1, \dots, G_n$ are finite
groups and $G= G_1*\dots* G_n$.  We consider various applications of
Theorem \ref{Tmain2}, beginning with augmentation ideals.  Here
$\bigoplus_{i=1}^n \Delta G_i \otimes_{\mathbb Z{G}_i}\mathbb Z{G} =
\Delta G$ (e.g.~\cite[\S 8.6]{Gruenberg70}), each $\Delta G_i$ is a
Swan lattice and $d_{G_i}(\mathbb Q \Delta G_i) =1$ (which makes (g2)
a triviality).  Thus Theorem \ref{Tmain2} gives
\begin{equation}\label{Efiniteaugmentation}
d_G(\Delta G) = \max_p \{ \sum_{i=1}^n
d_{G_i}(\Delta G_i/p\Delta G_i) \}.
\end{equation}
Hence \emph{the augmentation ideal of a free product of finite
groups is a Swan module}.  When $G_1, \dots, G_n$ have mutually
coprime orders, then \eqref{Efiniteaugmentation} with Lemma
\ref{Lmodulegenerators} gives Theorem \ref{Tfinitegroups}(a).

Next consider relation modules.  Let $R_i \hookrightarrow F_i
\twoheadrightarrow G_i$ be a free presentation and $R
\hookrightarrow F \twoheadrightarrow G$ the resulting free
presentation of $G$.  Then $\overline R = \bigoplus_{i=1}^n
\overline R_i \otimes_{\mathbb Z G_i}\mathbb Z G$ and
$d_{G_i}(\mathbb{Q} \overline R_i) = d(F_i)$
(see e.g.~\cite[Theorem 2.7]{Gruenberg76}).  From Theorem
\ref{Tmain2} we obtain

\begin{equation}\label{Efiniterelation}
d_G(\overline R) = \max_p\{ \sum_{i=1}^n
d_{G_i}(\overline R_i/p\overline R_i) \}.
\end{equation}
This has Theorem \ref{Tfinitegroups}(b) as corollary.  Furthermore,
the Gru\v sko-Neumann theorem tells us that a free
presentation of the free product of the groups $G_i$ is of the form
$*_i F_i \twoheadrightarrow *_i G_i$ with each $F_i$ a free group
mapping onto $G_i$, consequently \eqref{Efiniterelation} shows
that \emph{the relation module of a finite free presentation
of a free product of finite groups is a Swan module}.

The augmentation ideal of a group $G$ may be viewed as the first
kernel of a $\mathbb{Z}G$-free resolution of $\mathbb{Z}$, while a
relation module is the second kernel of such a resolution.  We
indicate how Theorem \ref{Tfinitegroups} may be generalised to
arbitrary kernels.

Let $H$ be a finite group, let
\begin{equation} \label{Eresolution1}
\dots \longrightarrow
(\mathbb{Z}H)^{e_2} \overset{\phi_2}{\longrightarrow}
(\mathbb{Z}H)^{e_1} \overset{\phi_1}{\longrightarrow}
\Delta H\overset{\phi_0}{\longrightarrow} 0
\end{equation}
be a $\mathbb{Z}H$-free resolution of $\Delta H$ and write
$K_i = \ker \phi_i$.  We recall a criterion for a lattice to be
a Swan module.  This is originally due to Swan, but
we use the improved form due to Jacobinski \cite{Jacobinski75}.

Let $\mathbb{Z}_{(H)}$ denote the semilocal subring $\{a/b \mid a,b
\in \mathbb{Z},\ (b,|H|) = 1\}$ of $\mathbb{Q}$ and take a minimal
free presentation $K \hookrightarrow (\mathbb{Z}_{(H)}H)^d
\twoheadrightarrow (K_s)_{(H)}$,
where $(K_s)_{(H)} = K_s \otimes \mathbb{Z}_{(H)}$.
Then $K_s$ is a Swan lattice if $\mathbb{Q}K$ contains a copy
of every non-trivial irreducible $\mathbb{Q}H$-module.  Since $H$ is a
finite group, every finitely generated $\mathbb{Q}H$-module is
uniquely a direct sum of irreducible $\mathbb{Q}H$-modules and we
have
\[
\mathbb{Q}K - (d - e_s + e_{s-1} - \dots + (-1)^se_1) \mathbb{Q}H +
(-1)^s \mathbb{Q} \Delta H = 0
\]
which shows that $\mathbb{Q} \Delta H$ is a summand of $\mathbb{Q}K$
whenever $s$ is odd and also when $s$ is even provided $d - e_s+
e_{s-1} - \dots + (-1)^se_1 \neq 1$; but if the sum is
1, then $K \cong \mathbb{Z}_{(H)}$ and $\mathbb{Z}$ has projective
period $\ell$ with $s+2 \equiv 0 \mod \ell$.

Finally $\mathbb{Q}K_s = (\mathbb{Q}H)^e \oplus (-1)^s
\mathbb{Q}\Delta H$, where
$e = e_s - e_{s-1} + \dots +(-1)^{s-1}e_1$.
Thus $\mathbb{Q} K_s$ needs $e+1$
$\mathbb{Q}H$-module generators if $s$ is even and $e$ generators
when $s$ is odd.  If $K_s$ fails to be a Swan lattice,
$\max_p\{ d_H(K_s/pK_s)\} = d = e+1 = d_H(\mathbb{Q}K_s)$.
Thus $d_H(K_s/pK_s)$ has constant value for all primes $p$.

\medskip

Now return to $G = G_1 * \dots * G_n$.  Let
\begin{equation} \label{Eresolution2}
\dots \longrightarrow
(\mathbb{Z}G_i)^{f_{i,2}} \overset{\theta_{i,2}}{\longrightarrow}
(\mathbb{Z}G_i)^{f_{i,1}} \overset{\theta_{i,1}}{\longrightarrow}
\Delta G_i\overset{\theta_{i,0}}{\longrightarrow} 0
\end{equation}
be a $G_i$-free resolution of $\Delta G_i$.  Applying
$\otimes_{\mathbb{Z}G_i} \mathbb{Z}G$ for each $i$, we obtain
a $\mathbb{Z}G$-free
resolution of $\Delta G_i \otimes_{\mathbb{Z}G_i} \mathbb{Z}G$.
The sum of these resolutions yields a $G$-free resolution
of $\Delta G$:
\begin{equation} \label{Eresolution3}
\dots \longrightarrow
(\mathbb{Z}G)^{f_2} \overset{\theta_2}{\longrightarrow}
(\mathbb{Z}G)^{f_1} \overset{\theta_1}{\longrightarrow}
\Delta G\overset{\theta_0}{\longrightarrow} 0,
\end{equation}
where $f_s = f_{1,s} + \dots + f_{n,s}$ and $\ker \theta_s =
\bigoplus_{i=1}^n \ker \theta_{i,s} \otimes_{\mathbb{Z}G_i}
\mathbb{Z}G$.

\begin{Prop} \label{Presolutionkernel}
Let $G = G_1 * \dots * G_n$ and take resolutions
\eqref{Eresolution2} and \eqref{Eresolution3}.  Suppose $s$
is a non-negative integer such that if $G_i$ has cohomological
period $\ell_i$, then $s+2 \not\equiv 0 \mod \ell_i$.  Then
\begin{enumerate}[\normalfont(a)]
\item
$\displaystyle
d_G(\ker \theta_s) = \max_p \{\sum_{i=1}^n d_{G_i}(\ker
\theta_{i,s}/p \ker \theta_{i,s})\}$;

\item
if $G_1, \dots, G_n$ have mutually coprime orders, then
\[
d_G(\ker \theta_s) = \max_k \{ d_{G_k}(\ker \theta_{k,s}) +
\sum_{i\ne k} d_{G_i}(\mathbb{Q}\ker \theta_{i,s})\}
\]
and $d_{G_i}(\mathbb{Q}\ker \theta_{i,s}) = f_{i,s} - f_{i,s-1}
+ \dots + (-1)^{s-1}f_{i,1} + \delta_s$, where $\delta_s = 1$ or
0 according as $s$ is even or odd.
\end{enumerate}
\end{Prop}
\begin{proof}
By our discussion of $H$ above, $\ker \theta_{i,s}$ is a Swan lattice
if $s$ is as stated.  Therefore Theorem \ref{Tmain2} yields (a); then
(b) is a consequence of (a), again using facts proved about $H$.
\end{proof}

Another application of Theorem \ref{Tmain2} in the context of finite
groups concerns finitely
generated nilpotent groups.  If $H$ is such a group, then there is a
prime $p$ so that the elementary abelian image $H/H'H^p$ has the
same minimal number of generators as does $H$.  Let us call such a
prime a \emph{generating prime} for $H$.
\begin{Prop} \label{Pgapnilpotent}
Let $H = H_1 * \dots * H_n$, where each $H_i$ is a finitely generated
nilpotent group.  If every $H_i$, except possibly one, is cyclic or
has $H_i/H_i'$ torsion-free, then $\gap(H) = 0$.
\end{Prop}
\begin{proof}
Suppose (without loss of generality) that $H_1$ is neither cyclic
nor has $H_1/H_1'$ torsion-free.  Let $q$ be a generating prime for
$H_1$ and $J(q)$ be the set of all $j$ such that $q$ is a generating
prime for $H_j$.  This absorbs all $H_j$ with $H_j/H_j'$
torsion-free.

Pick a generating prime $p_i$ for each $i$, but take $p_i = q$
whenever $i \in J(q)$.  Set $G_i = H_i/H_i'H^{p_i}$.  Then $G =
G_1* \dots *G_n$ is a homomorphic image of $H$ and $d(G) = d(H)$.
By reordering if necessary we may assume $J(q) = \{1,\dots,r\}$.
Then $G_i$ is cyclic for all $i > r$ and $p_i \not= q$.  Hence
$$\max_p\{ \sum_{i=1}^n d_{G_i}(\Delta G_i/p\Delta G_i) \} =
\sum_{i=1}^r d_{G_i}(\Delta G_i/q\Delta G_i) + n-r = d(G)$$
and so $\gap G = 0$ (by Theorem \ref{Tmain2}) whence $\gap H = 0$
by Lemma \ref{Lgap0}.
\end{proof}

It would be interesting to have a converse of Proposition
\ref{Pgapnilpotent}.  The closest we come in this direction is
Theorem \ref{Tinfinitegap} (cf.~Theorem \ref{Taugmentation}).

\section{Infinite Free Factors} \label{SIFF}

In order to use Theorem \ref{Tmain2} for infinite groups, we must find
good Swan modules.  We begin with augmentation modules.

\begin{Prop} \label{Paugmentation}
Let $A$ be a non-trivial finitely generated nilpotent group with
$A/A'$ torsion-free and let $G$ be a finite group such that
\begin{enumerate}[\normalfont(a)]
\item
$d(G) = d(G/G')$ and
\item
$d_G(\Delta G/p\Delta G) = d(G/G'G^p)$ for all $p\in \pi(G)$.
\end{enumerate}

If $H= G \times A$, then $\gap H = 0$, $\Delta H $ is both a Swan
lattice and a good $\mathbb Z H$-module for $\pi(H, \Delta H) =
\pi (G)$.  Furthermore $d(H) = d(G) + d(A)$ and $\delta_H(\Delta H)
= d(A)$.
\end{Prop}

For the proof of Proposition \ref{Paugmentation}, we require the
following result.

\begin{Lem}\label{Laugmentationmodp}
Let $G$ be a finite group, let $B$ be a finite nilpotent group and let
$E = G \times B$.  For any prime $p$,
\[
d_E(\Delta E/p\Delta E) = \max\{ d_G(\Delta G/p\Delta G), \dim G/G'G^p
+ \dim B/B'B^p\}.
\]
\end{Lem}
\begin{proof}
If $B(p)$ is the Sylow $p$-subgroup of $B$, $d_B(\Delta B/p\Delta B)
= d(B(p))$ (\cite[Lemma 7.17(i)]{Gruenberg76}) and $d(B(p)) =
d(B/B'B^p)$.  The result now follows from \cite[Theorem 2p]{CGK74}.
\end{proof}

Lemma \ref{Laugmentationmodp} shows that every finite nilpotent group
$G$ has property (b) of Proposition \ref{Paugmentation}.
There are plenty of
non-nilpotent finite groups having properties (a) and (b).  Let
$G= P \times B$, where $P=P'$, $B' = 1$, $\pi(P)\subseteq \pi(B)$ and
$d(B/B^p) \ge d_P(\Delta P)$ for all $p\in \pi(G)$.  For any
$k \geq 1$, the direct power $G^{(k)}$ of $G$ satisfies
$d_{G^{(k)}}(\Delta G^{(k)}/p\Delta G^{(k)}) = k\, d(B/B^p)$ for all
$p\in \pi(G)$ (again use \cite[Theorem 2p]{CGK74}).
Thus $G^{(k)}$ has
property (b) for all $k \geq 1$.  If $k$ is large enough, then by a
theorem of Wiegold (e.g.\ \cite[Theorem 6.13]{Gruenberg76}) $G^{(k)}
= kd(G/G') = d( G^{(k)} / G^{(k)'})$, whence $G^{(k)}$ has
property (a).

\begin{proof}[Proof of Proposition \ref{Paugmentation}]
When $G=1$ the result follows from Corollary \ref{Cgap0},
so henceforth we assume $G\ne 1$.
Set $\bar{A} = A/A^{|G|}$ and $\bar{H} = G \times \bar{A}$.
Now $\gap (G/G')= 0$ implies $\gap(G) = 0$, by \ref{Lgap0},
and this implies $\gap(\bar H )= 0$ (use \cite[(2.3)]{Gruenberg73}
and note that what is called in that paper the presentation rank
$\pr(E)$ of the finite group $E$ coincides with $\gap(E)$:
\cite[\S 12]{Gruenberg79}).  It will follow from Lemma \ref{Lgap0}
that $\gap(H)=0$ and that
$\Delta H$ is a Swan module, provided that $d(H) = d(\bar{H})$.

If $n\ge 2$ is a positive integer, $C$ is a finite abelian group such
that $nC=0$, and $D$ is a finitely generated free abelian group, then
$d(C\times D/nD) = d(C) + d(D)$.  Applying this with $n = |G|$, $C =
G/G'$ and $D = A/A'$, we find that
\[
d(H) \leq d(G)+d(A) = d(G/G')+d(A) = d\bigl((G/G') \times (\bar A/
\bar A')\bigr)\leq d(\bar H) \leq d(H).
\]
Thus $d(H) = d(\bar H)$ and also $d(H) = d(G) + d(A)$.

Consider condition (g1).  Let $a_1, \dots,a_e$ be free generators
of $A$.  By Corollary \ref{Cgap0},
$\{a_1-1, \dots, a_e-1\}$ is a minimal set
of $\mathbb{Z} A$-module generators of $\Delta A$ modulo $p\Delta A$
for every prime $p$.  Choose $p \not\in \pi(G)$.  Since $d_G(\Delta
G/p\Delta G) =1$ (by Lemma \ref{Lmodulegenerators}) there exists
$x \in \Delta G$ so that $x\mathbb Z G + p\Delta G = \Delta G$.  Let
\begin{equation}\label{Eaugmentation}
X = \{ x+(a_1-1)\hat G, \ a_2-1, \dots, a_e-1\}
\end{equation}
(where $\hat{G}$ is defined in the first paragraph of \S
\ref{SPreliminaries}).
We claim $X\mathbb Z H + p\Delta H= \Delta H$.  For if $z \mapsto
\bar{z} \colon \Delta H \twoheadrightarrow
\Delta H/p\Delta H$ denotes the natural epimorphism and $W =
\bar{X}\mathbb{Z}H$, then
\begin{equation} \label{Eaugmentation1}
(x + \hat{G}(a_1-1))(|G| - \hat{G}) = |G|x
\end{equation}
shows $\bar{x} \in W$, whence $\hat{G}(\bar{a}_1 - 1) \in
W$.  Also $\bar {g}-1 \in W$ for all $g \in G$ since
$\bar{x}\mathbb{Z}H$ contains $(\Delta G + p\Delta H)/p\Delta H$.
Hence $|G| - \hat{G} + p\Delta H \in W$, which gives
\[
|G|(\bar{a}_1 - 1) = (|G| - \hat{G} + \hat{G})(\bar{a}_1 - 1) \in W,
\]
and so $(\bar{a}_1 - 1) \in W$.  Thus $W = \Delta H/p\Delta H$.  We
conclude $d_H(\Delta H/p\Delta H) \leq e$ and hence equality because
$\Delta H \twoheadrightarrow \Delta A$.  Thus $\delta_H(\Delta H) = e$
and of course $d_H(\Delta H/p\Delta H) \geq e$ for \emph{every} prime
$p$, whence property (g1).

Finally property (g2).  Choose a finite set of primes $\pi$ properly
containing $\pi(G)$ and find a nested $\pi$-family $Y_1 \supseteq
Y_2 \supseteq \dots \supseteq Y_t$ in $\Delta G$ that satisfies (g2).
We may assume that there is a positive integer $s \le t$ such that
$p_i \in \pi(G)$ for $i<s$ and $p_i \notin \pi(G)$ for $i\ge s$.
Let $Y_i = \{ y_1, \dots, y_{d_i}\}$ (so $Y_t = \{y_1\}$).  Define
\[
X_t = \{ y_1+(a_1-1)\hat G, a_2-1, \dots, a_e-1\}
\]
for $i\ge s$, and for $i<s$ define
\[
X_i = \{y_1+(a_1-1)\hat G, a_2-1, \dots, a_e-1, a_1-1, y_2, \dots,
y_{d_i}\}.
\]
Then $X_1 \supseteq \dots \supseteq X_t$ and we claim this is a
nested $\pi$-family in $\Delta H$.  First consider $i<s$.
Since $X_i$ generates the same $\mathbb{Z}H$-submodule as $S :=
\{ y_1, y_2, \dots, y_{d_i}, a_1-1, a_2-1, \dots, a_e-1\}$
and clearly $S\mathbb Z H + p_i\Delta H = \Delta H$, we need to
show that $d_H(\Delta H/p_i\Delta H) = d_i + e$.  With $\bar{ H}
= G \times \bar {A}$ as above, $H \twoheadrightarrow \bar{ H}$
induces a bijection $S \isoto \bar {S}$ (since the kernel of $\Delta H
\twoheadrightarrow \Delta{\bar{ H}}$ is $(\Delta A^{|G|})
\mathbb{Z}H$).  So we are reduced to proving $d_{\bar{H}}
(\Delta{\bar{H}}/p_i \Delta {\bar{ H}}) = d_i + e$.  This is true
by property (b) on $G$ and Lemma \ref{Laugmentationmodp}.

Next consider the case $i \ge s$.
Since $X_i$ is the set $X$ of \eqref{Eaugmentation}
with $x = y_1$ and $p= p_i$, we see that $X_t$ minimally generates
$\Delta H$ modulo $p_i\Delta H$.  Simplify the notation now by setting
$X=X_i (=X_t)$ and $x=y_1$.  Let $q$ be a prime not in $\pi$.
Then (by (g2) for $\Delta G)$,
$d_G\bigl(\Delta G/(q\Delta G + x\mathbb{Z}G)\bigr)
\leq 1$, whence there exists $z \in \Delta G$ such that $\{x,z\}
\mathbb{Z}G + q \Delta G = \Delta G$.  Then $X\cup \{z\}$ generates
$\Delta H$ modulo $q\Delta H$.

To finish the proof, we must show that
$\Delta H/X \mathbb{Z}H$ has finite exponent.  Suppose $m$ is the
exponent of the finite module $\Delta G/x\mathbb ZG$.  Then $m$ is
prime to $p$ ($= p_i$) and \eqref{Eaugmentation1} shows that$|G|x
\in X\mathbb{Z}H$.  So $|G|m\Delta G$ is contained in
$X\mathbb{Z}H$, which gives $|G| m(|G| - \hat{G}) \in X\mathbb{Z}H$.
Also
\[
|G|\hat{G} (a_1-1)
= |G|(x + \hat G(a_1-1)) - |G|x \in X\mathbb{Z}H,
\]
whence $|G|^2m(a_1-1) = |G|m(|G|-\hat{G} + \hat{G})(a_1-1) \in
X\mathbb{Z}H$.  Therefore $|G|^2 m\Delta H \subseteq X\mathbb{Z}H$,
which completes the verification of (g2) and hence the proof of
Proposition \ref{Paugmentation}.
\end{proof}

The following result is an immediate consequence of Theorem
\ref{Tmain2} and Proposition \ref{Paugmentation}:
\begin{Thm} \label{Taugmentation}
Let $H = H_1* \dots * H_n$, where each $H_i$ is of the form $G_i
\times A_i$, $A_i$ is a finitely generated nilpotent group with
$A_i/A_i'$ torsion-free and $G_i$ is a
finite group such that
\begin{enumerate}[\normalfont(a)]
\item
$d(G_i) = d(G_i/G_i')$ and
\item
$d_{G_i}(\Delta G_i/p\Delta G_i) = d(G_i/G_i'G_i^p)$
for every $p\in \pi(G_i)$.  Then
\[
d_H(\Delta H) = \max_p\{\sum_{i=1}^n d_{H_i}
(\Delta H_i/p\Delta H_i)\}.
\]
\end{enumerate}
\end{Thm}

Theorem \ref{Tinfinitegap} is an immediate corollary of Theorem
\ref{Taugmentation} provided we note that $\delta_{H_i}
(\Delta H_i) = (1-\delta_{A_i,1})d(A_i) + \delta_{A_i.1}$.

\smallskip

We now turn to relation modules.  These are harder to deal with.  Let
$G$ be a non-trivial finite group, let $C$ be an infinite cyclic
group, and let $H = G \times C$.  Choose a finitely generated free
presentation $R \hookrightarrow F \twoheadrightarrow G$ for $G$ and
form the finitely generated free presentation $S \hookrightarrow
F*C \twoheadrightarrow H$ for $H$.
\begin{Prop} \label{Pg1}
$\overline{S}$ satisfies (g1) with $\delta_H(\overline{S}) = d(F)$.
\end{Prop}
\begin{proof}
Let $\alpha \colon F*C \twoheadrightarrow G*C$
denote the natural epimorphism induced by the
surjection of $F$ onto
$G$ and the identity on $C$, let $\beta \colon G*C
\twoheadrightarrow G \times C$ denote the
natural epimorphism induced by the identity on $G$ and $C$,
and let $\gamma \colon F*C \twoheadrightarrow G\times C$ denote
the natural epimorphism induced by the surjection of $F$ onto $G$
and the identity on $C$.  Set $R_0 = \ker \alpha$, $S= \ker\gamma$.
Then we have a commutative diagram with exact rows
\[
\begin{CD}
1 @>>> R_0 @>>> F*C @>\alpha>> G * C @>>> 1\\
@. @VVV @| @VV\beta V\\
1 @>>> S @>>> F*C @>\gamma>> G \times C @>>> 1
\end{CD}
\]
where the map $R_0 \to S$ is the natural inclusion.  Let
$E=\ker\beta$ and write $C = \langle c\rangle$.  By the snake
lemma (or the map $\beta$ can be identified with the natural
map $(F*C)/R_0 \to (F*C)/S$), we see that $S/R_0 \cong E$.  Also $E$
is the free group on $\{[g,c^i] \mid 1 \ne g \in G,\ 0 \ne i \in
\mathbb{Z}\}$, so $S = R_0 \rtimes E$, which gives the extension of
$\mathbb{Z}H$-modules
\begin{equation} \label{Eextension1}
R_0/[R_0,S] \rightarrowtail \overline{S} \twoheadrightarrow E/E'.
\end{equation}
Choose an inverse $t \colon G \to F$ to the free presentation $F
\twoheadrightarrow G$ (so $t(G)$ is a transversal for $R$ in $F$)
and assume that $t(1) = 1$.
Then $R_0 = *_w(w^{-1}Rw)$ (the free product of the subgroups
$w^{-1}Rw$) where $w$ runs through all those elements
of $F*C$ that can be expressed as a word in $t(G)$ and $C$,
and start with an element of $C$.  Hence
$\overline{R}_0 = \overline{R} \otimes_{\mathbb{Z}G}\mathbb{Z}[G*C]$
and then, because $R_0 \hookrightarrow S \twoheadrightarrow E$ is
exact, $R_0/[R_0,S] = (\overline{R} \otimes_{\mathbb{Z}G}
\mathbb{Z}[G*C]) \otimes_{\mathbb{Z}E} \mathbb{Z} \cong \overline{R}
\otimes_{\mathbb{Z}G} \mathbb{Z}H \cong \overline{R}
\otimes_{\mathbb{Z}} \mathbb{Z}C$; also $E/E' \isoto
\Delta G\otimes_{\mathbb{Z}} \mathbb{Z}C$ as $\mathbb{Z}H$-modules
via $[g,c^i] \mapsto (g-1) \otimes (c^i - 1)$ and $\Delta(C)
\isoto
\mathbb{Z}C$ via $c-1 \mapsto 1$.  Hence \eqref{Eextension1} becomes
\begin{equation} \label{Eextension2}
\overline{R} \otimes_{\mathbb{Z}} \mathbb{Z}C \hookrightarrow
\overline{S} \overset{\sigma}{\twoheadrightarrow} \Delta G
\otimes_{\mathbb{Z}} \mathbb{Z}C.
\end{equation}
Now $\tau \colon \Delta G
\otimes_{\mathbb{Z}} \mathbb{Z}C \to \overline{S}$ defined by $(g-1)
\otimes c^i \mapsto \overline{c^{-i}[t(g),c]c^i} =
\overline{[t(g),c]}c^i$ is a $\mathbb{Z}C$-homomorphism, whence $\tau
\hat{G}$ is a $\mathbb{Z}H$-homomorphism.
One checks that $(\tau \hat{G}) \sigma = |G| \id$.

So if $q$ is prime to $|G|$, then \eqref{Eextension2}
taken mod $q$ is a split extension.  The free presentation $R
\hookrightarrow F \twoheadrightarrow G$
determines the relation sequence $\overline{R} \rightarrowtail
\mathbb{Z}G^{d(F)} \twoheadrightarrow \Delta G$ and this splits
modulo $q$.  Thus
$\overline{R}/\overline{R}q \oplus \Delta G/q\Delta G
\cong (\mathbb{Z}G/q\mathbb{Z}G)^{d(F)}$ and therefore
\[
\overline{S}/q\overline{S} \cong (\mathbb{Z}G/q\mathbb{Z}G)^{d(F)}
\otimes_{\mathbb{Z}} \mathbb{Z}C/q\mathbb{Z}C \cong
(\mathbb{Z}H/q\mathbb{Z}H)^{d(F)},
\]
which gives $\delta_H(\overline S) = d(F)$.

Finally, let $p\in\pi(G)$ and write $E= F * C$.  Then $\overline S
\twoheadrightarrow SE'/E'$ is an $H$-homomorphism, $SE'/E'$ is
$H$-trivial and free abelian of rank $d(F)$, whence $\overline
S/p\overline S$ has an $H$-trivial $H$-image of
$\mathbb F_p$-dimension $d(F)$.  Thus $d_H(\overline
S/p\overline S) \geq d(F)$.
\end{proof}

Continuing with the notation of the last proof, take $C$-coinvariants
of \eqref{Eextension2}, giving
\begin{equation} \label{Eextension3}
\overline{R} \hookrightarrow \overline{S} \otimes_{\mathbb{Z}C}
\mathbb{Z} \twoheadrightarrow \Delta G,
\end{equation}
which is exact because $H_1(C,\mathbb{Z}C) = 0$.  The cocycle
determining the extension \eqref{Eextension2} is
$\psi \colon h \mapsto \tau h - \tau$.  If $h = gc^i$,
then $\psi(h) = \psi(g)$ because $\tau$ is a $C$-homomorphism.
The cocycle determining \eqref{Eextension3} is $\psi_C := \psi
\otimes_{\mathbb{Z}C} 1_{\mathbb{Z}}$.

Let $t(a)t(b) = t(ab) \rho(a,b)$ for $a,b \in G$.  Then
\[
\psi(g) \colon (a-1) \otimes 1 \mapsto \bigl( (ag^{-1} - g^{-1})
\otimes 1 \bigr) \tau g - \bigl( (a-1) \otimes 1 \bigr)\tau =
\bar{\rho}(a,g^{-1})g(1-c) \in \overline{S} \Delta C,
\]
whence $\psi_C(g) = 0$.  We conclude that \eqref{Eextension3} splits.

\begin{Lem} \label{LSwan}
Assume that $G$ is nilpotent and that there exists a relation prime of
$G$ that is also a $\Delta G$-prime.  Then $\overline{S}$ is a
Swan module and every $p$ that is both $\Delta G$-prime and
$\overline R$-prime is an $\overline{S}$-prime.
\end{Lem}
\begin{proof}
Let $p$ be an augmentation prime and also a relation prime for $G$.
Then $d_G(\Delta G) =
d(G/G'G^p)$ \cite[Lemma 7.17(i)]{Gruenberg76} and $d_G(\overline{R}) =
d(R/[R,F]R^p)$ \cite[Lemma 7.17(ii)]{Gruenberg76}.  Hence
\[
d_G(\overline{R}) + d_G(\Delta G) = d(R/[R,F]R^p) + d(G/G'G^p).
\]
Since \eqref{Eextension3} splits, there is an epimorphism
$(\overline{S} \otimes_{\mathbb{Z}C} \mathbb{Z})/p(\overline{S}
\otimes_{\mathbb{Z}C} \mathbb{Z}) \twoheadrightarrow
R/[R,F]R^p \oplus G/G'G^p$, which shows
\[
d(R/[R,F]R^p) + d(G/G'G^p) \le d_G(\overline{S}
\otimes_{\mathbb{Z}C}\mathbb{Z})/(\overline{S} \otimes_{\mathbb{Z}C}
\mathbb{Z})p \le d_H(\overline{S}/p\overline{S}) \le
d_H(\overline{S}).
\]
But $d_H(\overline{S}) \le d_G(\overline{R}) + d_G(\Delta G)$, by
\eqref{Eextension2}, so that we have achieved $d_H(\overline{S})
= d_H(\overline{S}/p\overline{S})$.
\end{proof}

Unfortunately the nilpotence of $G$ is not enough to ensure that
the remaining hypothesis in Lemma \ref{LSwan} is true:
cf.~\cite[Proposition 7.19]{Gruenberg76}.
But if $G$ is cyclic (as it is from
here on), then the set of all relation primes of $G$ coincides with
the set of all $\Delta G$-primes and this set is $\pi(G)$.
\begin{Lem} \label{Lg2}
Continue with the previous notation.  If $G$ is cyclic, then
$\overline{S}$ satisfies (g2) for $\pi(G\times C, \overline S) =
\pi(G)$.
\end{Lem}
\begin{proof}
Let $n = |G|$, write $G=\langle a\rangle$, let $F = \langle x\rangle$
(so $F$ is infinite cyclic), let $R = \langle x^n\rangle$, and let $x
\mapsto a$ in the presentation $R\hookrightarrow F \twoheadrightarrow
G$.  Then $S$ is the normal closure in $E = F*C$ of $\{x^n, [x,c]\}$,
whence $[S,E]$ is generated mod $E_3$
(the third term $[[E,E],E]$ of the lower central
series of $E$) by $[x^n,c] \equiv [x,c]^n \mod E_3$.  Thus
$[S,E] \not= E'$, since $E'/E_3$ is infinite cyclic on $[x,c]E_3$.
Therefore $d_H(\overline S) = 2$.  Hence $d_H(\overline{S}/
p\overline{S}) =2$ for every $p \in \pi(G)$ (by \ref{LSwan}).

The map $t\colon G \to F$ is here $t(a^i) = x^i$ for $0 \le i < n$
and $\tau$ maps $(a-1) \otimes 1$
to $\overline{[x,c]}$.  The calculation of $\psi(g)$ gives
\[
\tau \hat{G} \colon (a-1) \otimes 1 \mapsto
n\overline{[x,c]} + (\sum_{g\in G} \bar{\rho}(a,g^{-1})g)(1-c).
\]
The term $\sum_{g\in G} \bar{\rho}(a,g^{-1}) g = \overline{x^n}$.
Set $z = n\overline{[x,c]} - \overline{x^n}c = \bigl((a-1)
\otimes 1\bigr) \tau\hat{G} - \overline{x^n}$.  Since $\overline{S} =
\{\overline{x^n},\overline{[x,c]}\}\mathbb{Z}H$, we see that
$\overline{S} = \{z,\overline{[x,c]}\} \mathbb{Z}H$.  Now
\[
z(n -\hat{G}c) = n^2\overline{[x,c]} - n\overline{x^n}c -
\Bigl(\bigl( (a-1) \otimes 1 \bigr)\tau\hat{G} - \overline{x^n}\Bigr)
\hat{G}c = n^2\overline{[x,c]},
\]
which shows that $\overline{S} = z\mathbb{Z}H + q\overline{S}$ for any
$q$ prime to $n$; and
$\overline{S}/z\mathbb{Z}H$ has exponent dividing $n^2$.

If $\pi$ is a finite set of primes properly containing $\pi(G)$,
then $\{z, \overline{[x,c]}\} \supseteq \{z \}$ is a nested
$\pi$-family whose minimal element has the
required properties for (g2).
\end{proof}

Combining Proposition \ref{Pg1}, Lemma \ref{LSwan} and Lemma
\ref{Lg2} gives
\begin{Prop} \label{Pfinal}
If $F/R \isoto
G$ is a minimal free presentation of the finite cyclic
group $G$, $C$ is an infinite cyclic group and $(F*C)/S
\isoto G\times
C$, then $\overline{S}$ is a good Swan module for $G \times C$ with
$\pi(G\times C, \overline S) = \pi(G)$.

\end{Prop}
We call the free presentation in Proposition \ref{Pfinal} the
\emph{natural} minimal free presentation of $G \times C$.
Proposition \ref{Pfinal} and Theorem \ref{Tmain2} combine to give

\begin{Thm}\label{Tfinal}
Let $H = H_1 * \dots * H_n$, where each $H_i = G_i \times C_i$, $G_i$
is finite cyclic and $C_i$ is infinite cyclic.  Take natural minimal
free presentations $F_i/R_i \isoto H_i$ and let $F/R \isoto H$ be
the resulting free presentation of $H$.  Then
\[
d_H(\overline R) = \max_p\{\sum_{i=1}^n d_{Hi}
(\overline R_i/p\overline R_i)\}.
\]
\end{Thm}

Finally in this section, we discuss \cite[Proposition
2]{BridsonTweedale07} (which is a generalisation of \cite[Proposition
1]{BridsonTweedale07}).  Let $n$ be a positive integer, let
$\rho_n = \rho_n(x,t)$ be the word defined
by $\rho_n(x,t) = (txt^{-1})x(txt^{-1})^{-1}x^{-n-1}$, and let $Q_n =
\langle x,t \mid \rho_n = x^n = 1\rangle$.  Then \cite[Lemma
3]{BridsonTweedale07} shows that $Q_n$ is isomorphic to the
HNN-extension $(\mathbb{Z}/n\mathbb{Z} \times \mathbb{Z}/n\mathbb{Z})
*_{\phi}$, where $\phi$ maps the first factor isomorphically to the
second.  Let $R \hookrightarrow F
\twoheadrightarrow Q_n$ be the presentation determined by the
generators and relators for $Q_n$ from above.
Then, of course, the relation module $\overline{R}$ is generated as
a $\mathbb{Z}Q_n$-module by the images of $\rho_n,x^n$ in
$\overline{R}$.  Define $q_n = (n+1)^n - 1$ and $c_n = nq_n$.
Let $\bar{\rho}_n$ indicate the image of $\rho_n$ in $\overline{R}$.
Now \cite[Lemma 4]{BridsonTweedale07} proves that
$[(txt^{-1})^n,x^n] = x^{c_n}$ in $\langle x,t \mid \rho_n\rangle$
(here $[a,b]$ denotes the commutator $aba^{-1}b^{-1}$).
Furthermore since $[(txt^{-1})^n,x^n] \in R'$ and
$\overline{R}/\bar{\rho}_n\mathbb{Z}Q_n$ is generated as a
$\mathbb{Z}Q_n$-module by the image of $x^n$ in
$\overline{R}/\bar{\rho}_n\mathbb{Z}Q_n$, we see that
$\overline{R}/\bar{\rho}_n\mathbb{Z}Q_n$ is a group of exponent
dividing $q_n$.  We deduce that if $p$ is a prime that does not
divide $q_n$, then $\overline{R}=p\overline{R}+
\bar{\rho}\mathbb{Z}Q_n$.  We conclude that $\overline{R}$ is a good
$\mathbb{Z}Q_n$-module with $\delta_{Q_n}(\overline{R}) = 1$ and
$\pi(Q_n,\overline{R}) = \pi(q_n)$, where $\pi(q_n)$
indicates the set of primes dividing $q_n$.  Also if $n \ge 2$ and
$p$ is prime which divides $n$,
it is not difficult to see that $d_G(\overline{R}/p\overline{R})
= 2$ and hence $\overline{R}$ is a Swan module.
We can now apply Theorem
\ref{Tmain2} to obtain the following result.
\begin{Prop} \label{PBridson}
Suppose $r$ is a positive integer and $(q_{m_i},q_{m_j})
= 1$ for $1 \le i<j\le r$.  Assume that $m_i \ge 2$ for all $i$.
If $\overline{R}$ is the relation module for the group
$\Gamma : = Q_{m_1}* \dots *Q_{m_r}$
and presentation
\[
\langle x_{m_1}, t_{m_1}, \dots, x_{m_r},t_{m_r} \mid
\rho_{m_1}(x_{m_1},t_{m_1}),\dots,
\rho_{m_r}(x_{m_r},t_{m_r})\rangle,
\]
then $d_{\Gamma}(\overline{R}) = r+1$.
\end{Prop}
This result is essentially that of
\cite[Proposition 2]{BridsonTweedale07}, though there explicit
generators for the relation module are given.  Of course, one could
obtain further results by combining the previous paragraph with
Proposition \ref{Pfinal} and Theorem \ref{Tmain2}.

\section{Appendix} \label{SApp}

Throughout this section $G$ is a finite group.  Let $M$ be a
finitely generated $\mathbb{Z}G$-module.  As usual, $M_{(p)}$ denotes
the localization of $M$ at the prime $p$.  To establish that $M$
is a good $\mathbb{Z}G$-module with $\pi(G,M) = \pi(G)\cup \pi(M)$
we must prove

\begin{Thm}\label{TA}
Let $\pi$ be a finite set of primes properly containing $\pi(G,M)$.
Then there exists a nested $\pi$-family in $M$ whose minimal element
$X$ satisfies
\[
d_G\bigl( M/(X\mathbb{Z}G + pM)\bigr) \leq 1
\]
for every prime $p \notin \pi$.
\end{Thm}

We know that $M$ has property (g1) (by Lemma
\ref{Lmodulegenerators}) and $M/X\mathbb{Z}G$ is a finitely generated
abelian group with $p(M/X \mathbb{Z}G) = M/X \mathbb{Z}G$ for
$p \in \pi \setminus \{\pi(G)\cup \pi(M)\}$, whence $M/X \mathbb{Z}G$
is finite.  Thus Theorem \ref{TA} supplies what is missing of (g2).

We recall some relevant facts.  If $R$ is a subring of $\mathbb{Q}$
containing $1/|G|$, then $RG$ is a maximal order in $\mathbb{Q}G$,
every central idempotent of $\mathbb{Q}G$ belongs to $RG$ and every
$RG$-lattice is $RG$-projective \cite[respectively
Proposition 27.1, Theorem 26.20, Theorem 26.12]{CurtisReiner81}.
If $U$ is an irreducible
$\mathbb{Z}G$-module such that $pU=0$ for some $p \in \pi(G)'$,
then $U$ can be viewed as a $\mathbb{Z}_{(p)}G$-module, whence there
exists a unique primitive central idempotent $e$ such that $Ue\ne 0$.

\begin{Lem}\label{Lisolattices}
If $V,W$ are projective $\mathbb{Z}_{(p)}G$-lattices such that
$\mathbb{Q}V \cong \mathbb{Q}W$, then $V \cong W$
\cite[Theorem 32.1]{CurtisReiner81}.
\end{Lem}

\begin{Lem}\label{Lorder}
Let $M$ be a finite $\mathbb{Z}G$-module such that $M_{(p)} = 0$
for all $p \in \pi(G)$.  Then the composition factors of $M$ can
be arranged in any given order.
\end{Lem}

This follows from \cite[Lemma 9.3]{Swan70}.  Here is another proof,
shown to us by Steve Donkin.

\begin{proof}
It suffices to prove the result when $M$ has prime-power order, so
assume that $|M| = p^n$.  Then $p\notin \pi(G)$ and $M$ is a module
over $\mathbb{Z}_{(p)}G$.

Let $M = M_0 \supset M_1\supset \dots \supset M_d = 0$ be a
composition series and let $(U_1, \dots, U_d)$ be a (possibly new)
sequence of the composition factors.  We show by an induction on $d$
that there exists a composition series of $M$ producing factors
in this order.

When $d=1$ there is nothing to prove.  Suppose $d=2$.  If $pM= 0$,
then $M$ is semi-simple and so we are done.  If $pM \ne 0$, then
$pM$ must be simple and so is $M/pM$.  Then $M \twoheadrightarrow
pM$ has kernel $pM$, whence $M/pM \cong pM$ and again we are done.

Let $d > 2$ and $M/M_1 \cong U_r$.  If $r=1$, then induction gives a
new series on $M_1$ that produces $(U_2,\dots, U_d)$, completing the
proof.  Suppose $r>1$ and $U_1 \cong M_k/M_{k+1}$.  Then $k>0$ and by
induction we can find a composition series of $M_1$ with $U_1$ at
the top.  By the case $d=2$, there exists a composition series
$M=M_0\supset M_1' \supset M_2'$ with factors $(U_1, U_r)$.  Finally,
use induction on $M_1'$ to obtain a composition series of $M_1'$
with factors $(U_2, \dots, U_d)$.
\end{proof}

\begin{Lem}\label{Lauto}
Let $M$ be a finitely generated $\mathbb{Z}G$-module, let $p$ be a
prime not in $\pi(G,M)$ and let $U$ be an irreducible
$\mathbb{Z}G$-module such that $pU=0$.  Let $e$ be the primitive
central idempotent of $\mathbb{Q}G$ such that $Ue=U$ and assume
that $(\mathbb{Q}M)e$ is not irreducible.  If $\theta, \phi\colon
M\twoheadrightarrow U$ are $\mathbb{Z}G$-epimorphisms and $q$ is
an integer prime to $p$, then there exists a
$\mathbb{Z}G$-automorphism $\mu$ of $M$ such that $\theta = \phi\mu$
and $\mu \equiv \id \mod qM$.
\end{Lem}

\begin{proof}
Let $M'$ be the torsion submodule of $M$.  Since $p\nmid |M'|$,
we see that
$M_{(p)}$ is a $\mathbb Z_{(p)}G$-lattice and $M_{(p)}e \not= 0$.
If $(\mathbb{Q}M)e = V_1 \oplus V_2$ with $V_i \ne 0$, set
$M_i = V_i \cap M_{(p)}e$ and then
$M_1 \oplus M_2$ is a $\mathbb{Z}_{(p)}G$-submodule
of finite index in $M_{(p)}e$.  Therefore $M_{(p)}e = M_1
\oplus M_2$ by Lemma \ref{Lisolattices} and there exist epimorphisms
$\sigma_1\colon M_1 \twoheadrightarrow U$, $\sigma_2\colon M_2
\twoheadrightarrow U$.

Our proof follows that of \cite[Proposition 9.5]{Swan70} but with
some refinements.  Given $\mathbb{Z}_{(p)}G$-maps $\alpha_1,
\alpha_2$ of $M_1, M_2$ to $U$, we shall write $\alpha_1 \oplus
\alpha_2$ for the map $(x_1, x_2) \mapsto \alpha_1(x_1) +
\alpha_2(x_2)$.  Let $\theta_i$ be the map on $M_i$ induced by
$\theta$.  Clearly one of $\theta_1, \theta_2$ is an epimorphism.
Let us suppose it is $\theta_2$.

Since $M_1$ is $\mathbb{Z}_{(p)}G$-projective, there exists
$\psi\colon M_1 \to M_2$ so that $\theta_2\psi = \theta_1 - \sigma_1$.
Define $\tilde \psi \colon M_1 \oplus M_2 \to M_1 \oplus M_2$ by
$(x_1, x_2) \mapsto (0, \psi(x_1))$.  Then $\tilde \psi^2 =0$.
Since $M_{(p)}e$ is a direct summand of $M_{(p)}$, any endomorphism of
$M_{(p)}e$ can be considered as one of $M_{(p)}$ by defining it to
be zero on $M_{(p)} (1-e)$.  Then $\End_{\mathbb Z_{(p)}G}(M_{(p)})
\cong \mathbb{Z}_{(p)} \otimes \End_{\mathbb{Z}G}(M)$ shows
$\tilde \psi = (1/r) \otimes \tau$ for a suitable integer $r$ prime
to $p$ and $\tau \in \End_{\mathbb{Z}G}(M)$.

As $|M'|, r, q$ are all prime to $p$, we may choose $s$ so that
$qr|M'|s \equiv 1 \mod p$; then $t := |M'|s$ satisfies $tM'=0$.
Now $tM \cap M' =0$ shows $tM$ is a $\mathbb{Z}G$-lattice and $\tau$
restricts to an endomorphism of $tM$.
Since $(tM)_{(p)} \isoto M_{(p)}$ and $\tilde \psi^2 =0$,
we see that $\tau_{(p)}^2 = 0$, which gives
$\tau^2(tM) = 0$.  Define $\alpha = \id_M + (qt)\tau$.  This is a
$\mathbb{Z}G$-automorphism of $M$ and we readily check that
\[
(\sigma_1 \oplus \theta_2) \alpha_{(p)} = \theta_1 \oplus \theta_2.
\]

A similar procedure leads to $\psi' \colon M_2 \to M_1$ such that
$\sigma_1\psi' = \sigma_2 - \theta_2$, $\tilde \psi' \colon
(x_1, x_2) \mapsto (\psi '(x_2), 0)$, $r' \tilde\psi' = \tau '_{(p)}$
for some $r'$ coprime to $p$ and some $\tau' \in \End_{\mathbb{Z}G}M$.
Then the $\mathbb{Z}G$-automorphism $\alpha' = \id_M + (qt)\tau'$ of
$M$ satisfies
\[
(\sigma_1 \oplus \theta_2) \alpha '_{(p)} = \sigma_1 \oplus \sigma_2.
\]

We conclude that the $\mathbb{Z}G$-automorphism $\lambda =
\alpha^{-1}\alpha'$ of $M$ satisfies $(\theta_1 \oplus
\theta_2)\lambda_{(p)} = \sigma_1 \oplus \sigma_2$.  Similarly we can
find an automorphism $\nu$ such that $(\phi_1 \oplus \phi_2)\nu_{(p)}
= \sigma_1 \oplus \sigma_2$.  Then $\mu = \nu \lambda^{-1}$ is an
automorphism as required: for $\phi_{(p)}\mu_{(p)} = \theta_{(p)}$
implies $\phi\mu \equiv \theta \mod M'$ (because $(p, |M'|)
=1$) whence $\phi\mu = \theta$.
\end{proof}

\begin{Lem}\label{Llattice}
Let $L$ be a submodule of the $\mathbb{Z}G$-lattice $M$ such that
$M/L$ is finite of order prime to $|G|$.  Then there exists a
submodule $N$ of $M$ such that $M/N$ is a finite cyclic
$\mathbb{Z}G$-module with the same composition factors as $M/L$.
\end{Lem}

\begin{proof}
Since $M/L$ is the direct sum of its primary parts, we may assume
that $M/L$ has $p$-power order for some prime $p$,
so we may view $M/L$
as a $\mathbb Z_{(p)}G$-module.  Let $V_1,\dots, V_n$ be the
irreducible $\mathbb{Q}G$-modules which appear in $\mathbb{Q}M$.
Then we may choose a $\mathbb{Q}G$-submodule $V$ of $\mathbb{Q}M$
which is isomorphic to $\bigoplus_{i=1}^n V_i$.  Set $A= V \cap
M_{(p)}$.  Since $p$ is prime to $|G|$,
all $\mathbb{Z}_{(p)}G$-lattices
are projective.  Therefore we may write $M_{(p)} = A \oplus B$
for a suitable $\mathbb{Z}_{(p)}G$-submodule $B$ of $\mathbb{Q}M$,
because $M_{(p)}/A$ is a $\mathbb{Z}_{(p)}G$-lattice.
Also we see from
Lemma \ref{Lisolattices} that every submodule of $A$ is cyclic, and
if $A_0$ is a $\mathbb{Z}_{(p)}$-submodule of $A$ such that $A/A_0$
is finite, then every composition factor of $M/L$ is a
$\mathbb{Z}G$-homomorphic image of $A$.  Thus by induction on the
length of the composition series of $M/L$, we see that there is a
$\mathbb{Z}G$-submodule $C$ of $A$ such that $A/C$ has the same
composition factors as $M/L$.  Clearly $A/C$ is a cyclic
$\mathbb{Z}G$-module.  Now set $N = (C \oplus B) \cap M$.
Then $M/N \cong A/C$ as $\mathbb{Z}G$-modules and the result
follows.
\end{proof}

\begin{Lem}\label{Lfinal}
Let $L, N$ be submodules of the finitely generated
$\mathbb{Z}G$-module $M$ such that $M/L$, $M/N$ are
finite with the same composition factors and $\pi(M/L) \cap
\pi(G,M) = \emptyset$.  Assume that if $U$ is a composition factor
of $M/L$ with $pU=0$ and attached to the primitive central idempotent
$e$, then $(\mathbb{Q}M)e$ is not irreducible as $\mathbb{Q}G$-module.
Let $q$ be any integer prime to $|M/L|$.  Then there exists an
isomorphism $\theta\colon L \isoto N$ such that $\theta(x) \equiv x
\mod qM$ for all $x \in L$.
\end{Lem}
\begin{proof}
By Lemma \ref{Lorder} we may arrange the composition factors of
$M/L$ and $M/N$ so that they appear in the same order; in other words,
there exist chains of submodules $M = L_0 \supset L_1 \supset \dots
\supset L_d = L$ and $M = N_0 \supset N_1 \supset \dots \supset
N_d = N$ for some $d \geq 1$ such that, for each $i =0, 1, \dots,
d-1$, $L_i/L_{i+1}$ is $\mathbb{Z}G$-irreducible and $L_i/L_{i+1}
\cong N_i/N_{i+1}$.  We prove the result by induction on $d$.

By Lemma \ref{Lauto} there exists a $\mathbb{Z}G$-automorphism $\mu$
of $M$ such that $\mu$ induces $L_1 \isoto N_1$ and $\mu(x) \equiv x
\mod qM$ for all $x \in M$.  By induction there exists an isomorphism
$\nu\colon \mu L \isoto N$ such that $\nu(y) \equiv y \mod qN_1$ for
all $y \in \mu L$.  Then $\nu\mu\colon L\isoto N$ and $\nu\mu(x)
\equiv x \mod qM$ for all $x \in L$.
\end{proof}

\begin{proof}[Proof of Theorem \ref{TA}]
Choose a nested $\pi$-family ${\mathcal X}$ in $M$ and let $X$ be the
minimal element.  Let $e_1,\dots, e_n$ be the primitive central
idempotents in $\mathbb QG$ and define
\begin{align*}
I &= \{i \mid (\mathbb{Q}M)e_i \text{ is not irreducible}\},\\
J &= \{i \mid (\mathbb{Q}M)e_i \text{ is irreducible}\}.
\end{align*}
Set $e = \sum_{i \in I} e_i$ and let $M_1/X\mathbb ZG$ be the
$\pi'$-primary component of $M/X\mathbb{Z}G$.  Define the
$\mathbb{Z}G$-submodule $M_2$ of $M_1$ containing $X\mathbb{Z}G$
by $M_2/X\mathbb{Z}G = (M_1/X\mathbb{Z}G)e$.  Then each composition
factor of $M_2/X\mathbb{Z}G$ is attached to a unique $e_i$ with
$i \in I$.  The kernel of $M_2 \twoheadrightarrow M_2/X\mathbb{Z}G$
contains $M_2'$, the torsion group of $M_2$, because $\pi' \cap
\pi(G) = \emptyset$.

By Lemma \ref{Llattice} there exists $N \supseteq M_2'$ such that
$M_2/N$ has the same composition factors as $M_2/X\mathbb{Z}G$ and
$M_2/N$ is a cyclic $\mathbb{Z}G$-module.  Let $q$ be an integer
prime to $|M_2/X\mathbb{Z}G|$.  By Lemma \ref{Lfinal} there
exists an isomorphism $\theta\colon X\mathbb ZG \isoto N$ such
that $x \equiv \theta(x) \mod qM$ for all $x\in X$.  Thus
\[
X \equiv \theta(X) \mod pM
\]
for all $p \in \pi$.  So if we replace the subset $X$ of each
member of the nested $\pi$-family $\mathcal X$ by $\theta(X) =:Y$,
then we obtain a new nested $\pi$-family whose minimal element
is $Y$.  Observe $N = Y\mathbb{Z}G$ and $M/Y\mathbb{Z}G$ is finite
because $M/M_2$ is finite.  Also $(M_1/M_2)e=0$, which makes
$M_2/Y\mathbb{Z}G = (M_1/Y\mathbb{Z}G)e$.

We claim $M_1/Y\mathbb{Z}G$ is cyclic.  Since we already know
$(M_1/Y\mathbb{Z}G)e$ is cyclic, it remains to check
$(M_1/Y\mathbb{Z}G)(1-e)$ is also cyclic.
It suffices to do this for each primary component.  If $p$ occurs
in the primary decomposition, then $(M_1)_{(p)}(1-e) =
\sum_{i\in J} (M_1)_{(p)}e_i$, where $(\mathbb{Q}M)e_i$ for
$i\in J$ are distinct irreducible modules, whence
$(M_1)_{(p)}(1-e)$ is cyclic.  Therefore so is
$(M_1/Y\mathbb{Z}G)(1-e)$.  Finally, if $p \in \pi'$, then
$M/(pM + Y\mathbb{Z}G) = M_1/(pM + Y\mathbb{Z}G)$,
which gives our result.
\end{proof}

\bibliographystyle{plain}

\end{document}